\documentclass[preprint]{imsart}

\RequirePackage{bm}
\RequirePackage[OT1]{fontenc}
\RequirePackage{amsfonts}
\RequirePackage[numbers]{natbib}
\RequirePackage[colorlinks,citecolor=blue,urlcolor=blue]{hyperref}
\RequirePackage[centertags, leqno, sumlimits, nointlimits, namelimits]{mathtools}
\RequirePackage[thmmarks,thref,amsthm,amsmath]{ntheorem}
\RequirePackage{paralist}
\RequirePackage[english, noabbrev]{cleveref}
\RequirePackage{color}
\RequirePackage{graphicx}
\RequirePackage{epstopdf}
\RequirePackage{microtype}
\RequirePackage{bm}
\RequirePackage{float}
\RequirePackage{accents}
\startlocaldefs
\numberwithin{equation}{section}
\theoremstyle{plain}
\newtheorem{theorem}{Theorem}[section]

\newtheorem{corollary}[theorem]{Corollary}

\theoremstyle{remark}

\newtheorem{assump}{Assumption}
\newtheorem{remark}[theorem]{Remark}

\newtheorem{example}[theorem]{Example}

\endlocaldefs

\newcommand{\R}{\mathbb{R}} %Real numbers
\newcommand{\E}{\mbox{I\negthinspace E}} %Expectation
\newcommand{\lebesgue}{\ensuremath{\lambda\!\!\!\;\!\lambda}} % Lebesgue measure
\begin{document}
	
	\begin{frontmatter}
		\title{Dependence correction of multiple tests with applications to sparsity}
		\runtitle{Dependence correction of multiple tests}
		
		\begin{aug}
			\author{\fnms{Marc} \snm{Ditzhaus}$^1$}
			\and
			\author{\fnms{Arnold} \snm{Janssen}$^2$}\\[0.5em]
			{\footnotesize{1. Technical University of Dortmund, Germany\\ 
			2. Heinrich-Heine University Duesseldorf, Germany} }
			\footnotetext{Corresponding author: Marc Ditzhaus, Technical University Dortmund, Vogelpothsweg 87, 44221 Dortmund, Germany, email: 
				marc.ditzhaus@tu-dortmund.de}
		\end{aug}
		
		\begin{abstract}
			The present paper establishes new multiple procedures for simultaneous testing of a large number of hypotheses under dependence. Special attention is devoted to experiments with rare false hypotheses. This sparsity assumption is typically for various genome studies when a portion of remarkable genes should be detected. The aim is to derive tests which control the false discovery rate (FDR) always at finite sample size. The procedures are compared for the set up of dependent and independent $p$-values. It turns out that the FDR bounds differ by a dependency factor which can be used as a correction quantity. We offer sparsity modifications and improved dependence tests which generalize the Benjamini-Yekutieli test and adaptive tests in the sense of Storey. As a byproduct, an early stopped test is presented in order to bound the number of rejections. The new procedures perform well for real genome data examples.
		\end{abstract}
		
		
		
	\end{frontmatter}

	\section{Introduction and motivation}\label{sec:introduction}
	
	In genomics, but also in various other fields, several tests are applied simultaneously. Throughout, let $(H_i,p_i)_{1\leq i \leq m}$ denote a family of hypotheses $H_i$ for big data experiments with associated $p$-values $p_i$ and vector $p=(p_1,\ldots,p_m)$. If nothing else is said, let $p_i$ be uniformly distributed on $(0,1)$ for each $i\in I_0$. Here, $I_0$ is the set of indices belonging to true hypotheses, i.e., $i\in I_0$ iff $H_i$ is true. Analogously, introduce the index set $I_1=\{1,\ldots,m\}\setminus I_0$ of false hypotheses. Since the type $1$ error increases naturally when $m$ grows, corrections for multiple testing are needed. The classical Bonferroni correction, where the test size $\alpha$ is divided by the number $m$ of tests, is too conservative, especially for large $m$. \citet{Benjamini_Hochberg_1995} promoted to use the false discovery rate (FDR) as decision criteria for multiple testing procedures. The FDR is defined by  
	\begin{align*}
	\text{FDR}= \E\Bigl( \frac{V}{R} \Bigr)\quad\text{with the convention }\frac00 = 0,
	\end{align*}
	where $V$ is the number of false rejections and $R$ equals the number of all rejections. The step-up (SU) multiple testing procedure of \citet{Benjamini_Hochberg_1995}, in short denoted by BH procedure, is based on the linear critical values $\alpha_{j:m}=(\alpha j)/m$ and  rejects all hypotheses $H_i$ with
	\begin{align}\label{eqn:BH_rej_rule}
	p_i \leq \alpha_{R:m}, \quad \text{where } R = \max\{j: p_{j:m} \leq  \alpha_{j:m}\}
	\end{align}
	and $p_{1:m}\leq p_{2:m}\leq \ldots\leq p_{m:m}$ denote the order statistics. The FDR was first established under the following assumption:
	\\
	
	\textbf{Basic independence assumption (BI):} \textit{The uniformly distributed $p$-values $(p_i)_{i\in I_0}$ corresponding to true hypotheses are mutually independent as well as independent of the $p$-values $(p_i)_{i\in I_1}$ belonging to false hypotheses.}\\
	
	 At least under BI, the BH test controls the FDR for a pre-specified level $\alpha$. In fact, we have, see \citep{Benjamini_Hochberg_1995,FinnerRoters2001}, that
	\begin{align}\label{eqn:FDR_BH}
	\text{FDR}  = \frac{m_0}{m}\alpha,
	\end{align}
	where $m_0$ is the number of true hypotheses and, thus, equals the cardinality of $I_0$. The BH procedure was extended and modified several times in the literature \cite{BenjaminiETAL2006,BenjaminiYekutieli2001,HeesenJanssen2016,HellerETAL2009,Sarkar2008,Storey2002,StoreyETAL2004,StoreyTibshirani2003}. Moreover, it was shown to be still FDR-$\alpha$-controlling under \textit{positive regression dependence on each of a subset} (PRDS) \cite{BenjaminiYekutieli2001}. However, test statistics for gene expression data may be negatively correlated like multinomial distributions when gene material is partitioned.  As this example illustrates, specific dependence structure may be difficult to justify. \citet{GuoRao2008} stated: "\textit{It is almost impossible to check from biological principles or real data sets whether the underlying test statistics satisfy the assumption of independence or positive dependence of some type}".  This is in line with the statement of \citet{Efron2009}: "\textit{we expect the gene expressions to be correlated}", see also \citet{MoskvinaSchmidt2008}. Consequently, the question arises:
	\begin{itemize}
		\item How can the FDR be controlled under general dependent $p$-values?
	\end{itemize}
	\citet{BenjaminiYekutieli2001} pointed out that for the BH procedure 
	\begin{align}\label{eqn:FDR_BY}
	\text{FDR} \leq \frac{m_0}{m}\alpha \sum_{i=1}^{m} i^{-1}  
	\end{align}
	always holds, even under arbitrary dependence. Consequently, the modified and corrected BH procedure, labeled by BY, with 
	\begin{align}\label{eqn:alpha BY}
	\alpha_{j:m}^{\text{BY}} = \frac{j \alpha}{m \sum_{i=1}^{m}i^{-1}} 
	\end{align}
	always controls the FDR. However, this correction of logarithmic rate, $\log(m+1) \leq \sum_{i=1}^{m} i^{-1} \leq 1+ \log(m)$, leads to a substantial loss of power under dependence compared to the classical BH procedure, in particular, when $m$ is large. In this context we want to quote \citet{Owen2005} who remarked that the BY correction "\textit{could be very conservative}". For step-down tests with possibly less rejections \citet{GuoRao2008} provided an improvement of the correction factors by around 20-30$\%$. Unfortunately, the upper bound in \eqref{eqn:FDR_BY} cannot be improved in general for SU tests \cite{BenditkisETAL2018,GuoRao2008}. Although the examples attaining the upper bound have a constructive nature and are dealing with "\textit{extremely unusal distributions}" (as stated by \citet{BenjaminiHellerYekutieli2009}), there is no hope to get universal FDR-$\alpha$-controlling SU-tests. Throughout we will avoid special dependence assumptions. In contrast to that, we here focus on the case that we have mainly noise and just a few signals, a situation which is typically present in genome wide association studies.
	
	\subsection{Truncated BH procedures for rare signals}\label{sec:trunc_BH}
	Let us consider an experiment, e.g., a genome wide association study, with rare signals, i.e., there are only a few false hypotheses and $(m_0/m)$ is expected to be close to $1$. It seems plausible to anticipate that the fraction $R/m$ of rejections is small as well. Theoretically, this can be supported by an observation of \citet{DitzhausJanssen2019} under independence that $m_0/m\to 1$ implies $R/m \to 0$ for the BH procedure. When the fraction $R/m$ is small the majority of the large critical values does not have an impact. Having this in mind, we suggest to focus on a smaller pre-specified portion $1 \leq j \leq k$ with $k \ll m$ of the  critical values $\alpha_{j:m}$ and to reduce the remaining ones, e.g., by truncation: $\alpha_{j:m}= \alpha_{k:m}$ for $j>k$. Among other, we will show that for these truncated critical values the correction factor reduces to $\sum_{i=1}^{k} i^{-1}$ and, hence, the critical values 
	\begin{align}\label{eqn:def_beta}
	\beta_{j:m}^{(k)} = \frac{\alpha}{m} \frac{\min(j,k)}{\sum_{i=1}^{k} i^{-1} }
	\end{align}
	always lead to FDR control by $(m_0/m)\alpha$. Depending on the choice of $k$, this correction factor can be a significant improvement compared to the BY correction  factor $\sum_{i=1}^{m} i^{-1}$, see Table \ref{tab:corr_fact} for illustration. As for the classical Bonferroni test ($k=1$), the SU test given by \eqref{eqn:def_beta} may lead to $R>k$ rejections.
	\begin{table} 
		\large
		\centering
		\caption{The correction factors under sparsity for various choices of $k$}\label{tab:corr_fact}
		\begin{tabular}{ c | c | c | c | c | c }
			$k$ & $20$ & $50$ & $100$ & $200$ & $m=5\cdot 10^5$ \\
			\hline
			${\sum_{i=1}^{k} i^{-1}}/{\sum_{i=1}^{m} i^{-1}}$ & $0.262$ & $0.328$ & $0.378$ & $0.429$ & $1.000$\\
			\hline
			$\sum_{i=1}^{k} i^{-1}$ & $3.597$ & $4.499$ & $5.187$ & $5.878$ & $13.699$ \\
			\hline
			$\log(k+1)$ & $3.044$ & $3.931$ & $4.615$ & $5.303$ & $13.122$ \\
			%\hline
			%$\log\log k$ & $1.0972$ & $1.3641$ & $1.5272$ & $1.6674$ & $2.5743$ {\color{red} Rausnehmen!}\\
		\end{tabular}
	\end{table}

	\subsection{Overview}
	The present paper aims to compare the FDR under
	\begin{itemize}
		\item the basic independence model as a benchmark and
		
		\item the general dependence setting.
	\end{itemize}
	Our FDR bounds differ by a factor which can be used as a correction factor for overall valid dependence procedures, similarly to \eqref{eqn:FDR_BY} and \eqref{eqn:alpha BY}.
	
	The paper is organized as follows. In Section \ref{sec:towards_SU_sparsity}, we first give an overview of existing step-up procedures leading to FDR-$\alpha$-control under the basic independence assumption. Then we adapt the truncation idea from the previous section to these procedures and explain which correction factors are needed to preserve the FDR-$\alpha$-control under dependence. By these correction factors we can quantify the price to pay for dependence. In this context, we also discuss adaptive procedures like the one of \citet{Storey2002}. Moreover, we introduce an early stop procedure, for which the number of rejections is upper bounded by a pre-specified limit. Section \ref{sec:main} contains our main Theorems generalizing the results from Section \ref{sec:towards_SU_sparsity}. A general upper FDR bound for data dependent critical values is presented in Theorem \ref{theo:doehler_exte}. Special attention is devoted to new designed SU-tests for sparsity models. Their advantage is a reduction of the dependence correction factor. In this context, Example \ref{example:extreme_dependence} explains how Storey's adaptive test can be modified. In Section \ref{sec:real_data}, we illustrate the new procedures' benefits by discussing two real data examples. In particular, the role of the truncation parameter $k$, see \eqref{eqn:def_beta}, is illustrated. The real data example of \citet{needleman:etal:1979} suggests to design SU-tests with enlarged critical values $\alpha_{j:m}$ for very small $j$. To overcome this difficulty, Theorem \ref{theo:new_alpha} offers another promising overall FDR-$\alpha$-controlling modification of the Benjamini-Yekutieli procedure. All proofs are collected in Section \ref{sec:proofs}.
	
	\section{First sparsity SU tests under dependence}\label{sec:towards_SU_sparsity}
	\subsection{Procedures under the basic independence assumption}
	We first present some known as well as a new SU procedure being (asymptotic) FDR-$\alpha$-controlling. In the next section, we explain how to correct the corresponding critical values to preserve the FDR-$\alpha$-control under any dependence structure. \\
	The linear critical values introduced in Section \ref{sec:introduction} are replaced by general (ordered) critical values
	\begin{align}\label{eqn:alpha_<...}
	0=\alpha_{0:m} <  \alpha_{1:m} \leq \alpha_{2:m} \leq \ldots \leq \alpha_{m:m} 
	\end{align}
	with the restriction $\alpha_{m:m}<1$. Then the SU test procedure and the number $R$ of rejections is given by \eqref{eqn:BH_rej_rule}. The number $V$ of falsely rejected hypotheses is defined by the sum of  all rejections of true nulls
	\begin{align*}
	V = \sum_{i\in I_0} \mathbf{1}\{p_i\leq \alpha_{R:m}\}.
	\end{align*} 
	A lot of critical values may be introduced by a non-decreasing, left continuous generating function
	\begin{align}\label{eqn:def_g}
		g:[0,1]\to[0,\infty]\quad\text{ with } g(0)=0<g(x) \text{ for } x>0.
	\end{align}
	The function $g$ can be regarded as the inverse of so-called rejection curves \cite{FinnerETAL2009}. The appertaining deterministic critical values are then
	\begin{align}\label{eqn:def_alpha_g}
	\alpha_{j:m} = g\Bigl( \frac{j}{m} \Bigr),\quad 1\leq j \leq m \text{ if }0<g\Bigl( \frac1m  \Bigr)\;\text{ and }\;g(1)<1.
	\end{align}
	In the following we list some prominent examples of such critical values. Let $0<\lambda < 1$ be a fixed truncation or tuning parameter.
	\begin{enumerate}
		\renewcommand\labelenumi{\bfseries(\text{W}\theenumi)}
		\item\label{enu:proc_BH} The BH procedure \citet{Benjamini_Hochberg_1995} is given by
		\begin{align*}
		\alpha_{j:m} = \frac{j\alpha }{m}\quad \text{and}\quad g(x)=\alpha x.
		\end{align*}
		
		\item\label{enu:proc_AORC} The procedure of \citet{FinnerETAL2009}, based on the asymptotic optimal rejection curve (AORC), corresponds to
		\begin{align*}
		\alpha_{j:m} = \frac{ j\alpha  }{ m -j(1-\alpha)}\quad (j<m)\quad \text{and}\quad g(x)= \frac{\alpha x}{ 1 - x(1-\alpha)}.
		\end{align*}
		Since  we have $g(1)=1$ here, the last coefficient needs to be modified such that $\alpha_{m:m}\in(\alpha_{m-1:m},1)$, we refer to \citet{FinnerETAL2009} and \citet{Gontscharuk2010} for a detailed discussion.
		
		\item\label{enu:proc_Blan+Roq} The procedure of \citet{BlanchardRoquain2008,BlanchardRoquain2009} is defined by
		\begin{align*}
		\alpha_{j:m} = \Bigl( (1-\lambda) \frac{ j\alpha  }{ m + 1 -j}\Bigr)\wedge \lambda \quad \text{and}\quad g(x)= (1-\lambda) \frac{ \alpha x}{ 1 + \frac1m - x} \wedge \lambda.
		\end{align*}
		
		\item \label{enu:proc_mod_Blan+Roq} A new procedure is given by a combination of (W\ref{enu:proc_AORC}) and (W\ref{enu:proc_Blan+Roq}): 
		\begin{align*}
		\alpha_{j:m} = \Bigl( (1-\lambda) \frac{ j \alpha}{ m - j(1-\alpha)} \Bigr) \wedge \lambda,\quad g(x)= (1-\lambda) \frac{ \alpha x}{ 1  - x(1-\alpha)} \wedge \lambda.
		\end{align*}
	\end{enumerate}
	It is well-known under the basic independence model that the procedures (W\ref{enu:proc_BH}) and (W\ref{enu:proc_Blan+Roq}) are FDR-$\alpha$-controlling, for the latter see Theorem 9 of \citet{BlanchardRoquain2009}. The AORC procedure (W\ref{enu:proc_AORC}) does not control the FDR under (BI) because the first coefficient $\alpha_{1:m}> 1/m$ is already too large. However, these values can be modified to achieve asymptotic FDR-$\alpha$-control, see \citet{FinnerETAL2009} and see also \citet{HeesenJanssen2015} for some modifications. The new combination approach (W\ref{enu:proc_mod_Blan+Roq}) is FDR-$\alpha$-controlling as well:
	\begin{theorem}\label{theo:mod_Blan+Roq} 
		Let $\alpha/(m-1)<\lambda<1$. Suppose that $m\geq \frac{1-\alpha}{\alpha - \alpha^2/4}$ and $m\geq (1-\lambda)\lambda/\alpha$. Then the procedure (W\ref{enu:proc_mod_Blan+Roq}) is FDR-$\alpha$-controlling under the basic independence assumption. 
	\end{theorem}
	  Various FDR-$\alpha$-controlling extensions have been established for "positive \mbox{dependent"} $p$-values, more precisely PRDS, reverse martingale models or extended BI-models with stochastically larger $p$-values than the uniform ones for $i \in I_0$, i.e., $P(p_i \leq x) \leq x$. 
	\\
	Since the BH procedure does not exhaust the FDR level $\alpha$ completely, compare to \eqref{eqn:FDR_BH}, the level $\alpha$ in the linear BH critical values were replaced by a data depended, adjusted level $\alpha'= \alpha (m/\widehat m_0)$, where $\widehat m_0$ is an estimator for $m_0$. These so-called adaptive procedure are expected to exhaust the level better because, heuristically, $\text{FDR}_{ m}\approx (m_0/m)\alpha' \approx \alpha$. Various estimators are suggested in the literature  \cite{benjamini:hochberg:2000,BenjaminiETAL2006,BlanchardRoquain2008,BlanchardRoquain2009,HeesenJanssen2016,schweder:spjotvoll:1982,StoreyTibshirani2003,zeisel:ETAL:2011}. Exemplarily, we present the formula of the Storey estimator \cite{Storey2002,StoreyETAL2004}:
	\begin{align}\label{eqn:def_Storey_estimator}
	\widehat m_0^{\text{Stor}} = m\frac{ 1 - \widehat F_m(\lambda) + \frac{1}{m} }{ 1- \lambda},
	\end{align}
	where $\lambda \in(0,1)$ is a tuning parameter and $\widehat F_m$ denotes the empirical distribution function of the $p$-values. Plugging-in this estimator into the BH procedure leads to FDR-$\alpha$-control \cite{StoreyETAL2004}. Conditions for FDR-$\alpha$-control of linear adaptive tests were reviewed and established by \citet{HeesenJanssen2015,HeesenJanssen2016}. To include these and more general adaptive procedures, we consider, additionally to the deterministic critical values, the data driven critical values 
	\begin{align}\label{eqn_alpha_g_lambda}
	\widehat \alpha_{j:m} = \min \Bigl( g\Bigl( \frac{j}{\widehat m_0} \Bigr), \lambda \Bigr),
	\end{align}
	where $\widehat m_0$ is any estimator of $m_0$. Some of our results are even valid for arbitrary, data-driven critical values
	\begin{align}\label{eqn:hatalpha_<...}
	0=\widehat \alpha_{0:m} < \widehat \alpha_{1:m} \leq \widehat\alpha_{2:m} \leq \ldots \leq \widehat\alpha_{m:m}. 
	\end{align}
	
	\subsection{Sparsity modifications}\label{sec:spars_mod}
	We extend here the truncation idea introduced in Section \ref{sec:trunc_BH} from the classical BH procedure to all the other procedures mentioned in the previous section. Thus, let $k\in\{1,\ldots,m\}$ be pre-chosen and let $g$ be one of the generating functions corresponding to the procedures (W\ref{enu:proc_BH})--(W\ref{enu:proc_mod_Blan+Roq}). The main focus lies on the sparse signal case, where $k\ll m$ is chosen, but the results hold for any $k$. At the end of the previous section, we introduced adaptive procedures using an estimation step. In the sparse case, we expect $Cm\leq m_0$ with $C$ close to $1$. Throughout the remaining paper, we consider just estimators $\widehat m_0$ for $m_0$ such that
	\begin{align}\label{eqn:cond_hatmo}
	mC\leq \widehat m_0 \leq \frac{m}{\delta}\quad \text{ for some }C,\delta\in (0,1].
	\end{align}
	We want to point out that the Storey estimator \eqref{eqn:def_Storey_estimator} may become larger than $m$. That is why we consider also the case $\delta<1$, whereas the choice $\delta=1$ may be more plausible in applications. By
	\begin{align*}
	\widehat m_0 = \max( Cm, \min( m\delta^{-1}, \widetilde m_0))
	\end{align*} 
	we can transfer any estimator $\widetilde m_0$ into an estimator fulfilling \eqref{eqn:cond_hatmo} for any pre-chosen $C,\delta$. 
	In the deterministic case $\widehat m_0=m$, we set $\delta=C=1$. In the spirit of Section \ref{sec:trunc_BH}, we truncate the critical values  \eqref{eqn_alpha_g_lambda}:
	\begin{align}\label{eqn:def_alpha_k}
	\widehat \alpha_{j:m}^{(k)} = \min \Bigl( g\Bigl( \frac{j\wedge k}{\widehat m_0} \Bigr), \lambda \Bigr).
	\end{align}
	For $g(x)=\alpha x$ and $\widehat m_0=m$ we obtain the Bonferroni as well as the BH tests by setting $k=1$ and $k=m$, respectively. As in \eqref{eqn:def_beta}, the critical values need to be corrected under dependence.  As an immediate consequence of our main Theorem \ref{theo:main_res}, we get:
	\begin{corollary}\label{cor:motiv:beta}
		The corrected critical values 
		\begin{align}\label{eqn:beta_mot}
			\beta_{j:m}^{(k)} = D_k^{-1}C_k^{-1} \widehat \alpha_{j:m}^{(k)}
		\end{align}
		specified by the values in Tables \ref{tab:Ck} and \ref{tab:Dk2} always lead to FDR-$\alpha$-control.
	\end{corollary} 
	The correction in \eqref{eqn:beta_mot} can be separated into two parts: (1) a procedure correction $C_k$ listed in Table \ref{tab:Ck} and (2) a dependence structure correction $D_k$, which is also affected by the choice between deterministic and adaptive critical values, see Table \ref{tab:Dk2}. 
	\begin{example}\label{exam:} 
		Let $k\sim m^\gamma$ for large $m$ and $0<\gamma<1$. Let us just discuss the BH procedure under dependence. Then 
		\begin{align*}
		D_k = \sum_{i=1}^k\frac{1}{i} \sim 1 + \log(k) \sim  \gamma (1+\log m) + (1-\gamma) \sim \gamma \sum_{i=1}^m\frac{1}{i} + (1-\gamma).
		\end{align*}
		Consequently, by truncation we reduce the dependence correction by a factor close to $\gamma$ compared to the classical BY correction.
	\end{example}
	
	\begin{table} 
		\large \centering
		\caption{The procedure correction factor $C_k$ for (W\ref{enu:proc_BH}), (W\ref{enu:proc_AORC}), (W\ref{enu:proc_Blan+Roq}),  (W\ref{enu:proc_mod_Blan+Roq}) with $B=k/(Cm)$} \label{tab:Ck}
		\begin{tabular}{ c || c  c  c c }
			& (W\ref{enu:proc_BH}) & (W\ref{enu:proc_AORC}) & (W\ref{enu:proc_Blan+Roq})& (W\ref{enu:proc_mod_Blan+Roq}) \\[0.2cm]
			\hline \hline\\[-0.4cm]
			$C_k$ & 1 &  $\frac{1}{ 1 - B(1-\alpha)}$ &  $\frac{1-\lambda}{ 1 + \frac1m - B} \wedge \frac{1}{\alpha}$ & $\frac{1-\lambda}{ 1 - B(1-\alpha)} \wedge  \frac{1}{\alpha}$ 		
		\end{tabular} 
	\end{table}
	\begin{table} 
		\large \centering
			\caption{Comparison of correction factor $D_k$ for the different dependence models}\label{tab:Dk2}
		\begin{tabular}{ c || c  c    }
			$D_k$& BI & Dependence  \\[0.1cm]
			\hline \hline\\[-0.3cm]
			deterministic & 1 &   $\sum_{i=1}^k i^{-1}$ \\[0.1cm]
			\hline\\[-0.3cm]
			adaptive & $C^{-1}$& $C^{-1}\log ( 1 + k/(C\delta))$
		\end{tabular} 
	\end{table}

	\subsection{Early stopped multiple procedures} \label{sec:early_stopp}
	In many studies, e.g., dealing with genes, a screening step is necessary in a first experiment.  But what shall be done when in the first experiment too many hypotheses are rejected and resources for a detailed examination are limited? Say, $R$ exceeds $k$, whereas the close investigation is limited to $k$ genes, for example. Ad hoc, it seems plausible to select just the hypotheses corresponding the smallest $k$ $p$-values. But, already under BI, well-known examples show that this new reduced  procedure is not FDR-$\alpha$-controlling. Our truncation idea cannot solve this limitation problem since $R>k$ is possible. However, our early stopped procedure ensures $R\leq \kappa$ for any pre-specified $\kappa\in\{1,\ldots,m-1\}$, at least whenever the very conservative Bonferroni test with $\alpha_{j:m}=(\alpha/m)$ rejects not more than $\kappa$ hypotheses. The data dependent critical values of the new procedure are given by
	\begin{align}\label{eqn:def_alpha_early}
	\widehat \alpha_{j:m} = \frac{ \alpha}{m} \min\{j,j^*\}\quad \text{with }j^*= 1 \vee \max\{i:\frac{i\alpha}{m}< p_{\kappa+1:m} \}.
	\end{align}
	\begin{theorem}\label{theo:early_stopp}
		Consider the critical values \eqref{eqn:def_alpha_early}. 
		\begin{enumerate}[(a)]
			\item\label{enu:theo:early_stopp_BI} Under BI we have $\text{FDR}\leq (m_0/m)\alpha$.
			
			\item\label{enu:theo:early_stopp_dep} Under arbitrary dependence models the corrected critical values
			\begin{align*}
			\widehat\beta_{j:m}^{(k)} = \min\{\widehat \alpha_{j:m},\widehat \alpha_{k:m}\}/\sum_{i=1}^k \frac{1}{i}
			\end{align*}
			lead to an FDR-$\alpha$-controlling procedure.
		\end{enumerate}
	\end{theorem}

	\section{Main results}\label{sec:main}
	\subsection{Procedures based on general generating functions}
	Here, we extend the results of Section \ref{sec:early_stopp} to critical values of the shape \eqref{eqn:def_alpha_k} for general generating functions \eqref{eqn:def_g} and estimators $\widehat m_0$ fulfilling \eqref{eqn:cond_hatmo}. The tuning parameter $\lambda$ can be freely chosen, it is also possible to set $\lambda=1$ whenever $g(k/(mC))<1$. As an extension of the procedure correction $C_k$, let us introduce via the right continuous inverse $g^{-1}$ of $g$:
	\begin{align}\label{eqn:Ck_ad+det}
	C_k^{\text{ad}} = \frac{1}{\alpha} \sup_{\delta/m\leq t \leq B} \frac{ g(t)}{t} \quad \text{ and }C_k^{\text{det}} = \frac{1}{\alpha} \sup_{j=1,\ldots,\lfloor mB\rfloor} \frac{ g(j/m)}{j/m},
	\end{align}
	where $B= \min(g^{-1}(\lambda) , k/(Cm))$ and $\lfloor x \rfloor$ is the integer part of $x\in\R$. In comparison to Section \ref{sec:spars_mod}, the procedure correction also depends on the choice between deterministic and adaptive critical values. We want to point out that most of all generating functions $g$ are convex. For these, both suprema are attained at the end point of the considered interval.
	\begin{theorem}\label{theo:main_res}
		Consider the SU-tests with critical values \eqref{eqn:def_alpha_k}.
		\begin{enumerate}[(i)]
			\item \label{enu:theo:main_res_BI} (BI bound) Suppose that $\widehat m_0 = \widehat m_0( (\widehat F_m(t))_{t\geq \lambda})$ only depends on the $p$-values $p_i\geq \lambda$. Then
			\begin{align}\label{eqn:main_res_BI}
			\text{FDR} \leq \frac{m_0}{m} \alpha C^{-1} C_k^{\text{ad}}.
			\end{align}
			In the deterministic case $\widehat m_0=m$, the inequality \eqref{eqn:main_res_BI} is still valid when $C_k^{\text{ad}}$ is replaced by $C_k^{\text{det}}$ and $C$ is set equal to $1$.
			
			\item\label{enu:theo:main_res_dep_det} (Arbitrary dependence, deterministic case) Suppose that $\widehat m_0=m$ and, in particular, $C=1$. Then
			\begin{align*}
			\text{FDR} \leq \frac{m_0}{m} \alpha\Bigl( \sum_{1\leq i\leq mB} \frac{1}{i} \Bigr) C_k^{\text{det}}
			\end{align*}
			provided that $( \sum_{1\leq i\leq mB} 1/i ) C_k^{\text{det}}<1$.
			
			\item\label{enu:theo:main_res_dep_adap} (Arbitrary dependence and adaptivity) Suppose additionally that $g$ is absolutely continuous, i.e., there exists a Lebesgue almost sure derivative $g'$ such that $g(x) = \int_0^x g'(t) \,\mathrm{ d }t$. In contrast to \eqref{enu:theo:main_res_BI}, we allow that the estimator $\widehat m_0=\widehat m_0((\widehat F_m(t))_{0\leq t \leq 1})$ may use the information of all $p$-values. Then 
			\begin{align}\label{eqn:FDR_upper_log}
			&\text{FDR}\leq \alpha \frac{m_0}{m} C^{-1}\Bigl( 1 + \log \frac{mB}{\delta} \Bigr) C_k^{\text{ad}}\\
			&\text{whenever }\frac{m}{\delta}g\Bigl( \frac{\delta}{m} \Bigr) + \int_\delta^{mB} \frac{1}{z} g'\Bigl( \frac{z}{m} \Bigr)\,\mathrm{ d }z < 1.\nonumber 
			\end{align}
		\end{enumerate}
	\end{theorem}
	The preceding Theorem provides a construction principle for SU-tests controlling the FDR at level $\alpha$. For this purpose, let us return to the dependence structure correction factor $D_k$, see Table \ref{tab:Dk2}. This factor depends, as the procedure correction $C_k\in\{C_k^{\text{det}},C_k^{\text{ad}}\}$, on the choice between deterministic and adaptive critical values. Consequently, we get
	\begin{corollary}\label{cor:}
		The corrected critical values
		\begin{align}\label{eqn:tildealpha}
		\widehat \beta_{j:m}^{(k)}:= C_k^{-1} D_k^{-1} \widehat \alpha_{j:m}^{(k)} 
		\end{align}
		lead to FDR-$\alpha$-controlling procedures. 
	\end{corollary}
	Note that \eqref{eqn:tildealpha} is scale invariant regarding $g$ and Theorem \ref{theo:main_res} always applies.

	\subsection{General FDR upper bound under dependence}
	
	\citet{Doehler2018} derived general upper FDR bounds for deterministic $\alpha_{j:m}$. For our purposes, we need to, among others, extend these bounds to data dependent critical values $\eqref{eqn:hatalpha_<...}$. Suppose
	\begin{assump}\label{ass:constant_alpha}
		For each $i\in \{1,\ldots,m\}$ let $p_i\mapsto (\widehat \alpha_{j:m}(p))_{1\leq j\leq m}$ be constant on the set $[0,\widehat \alpha_{R:m}]$, i.e., the set on which $H_i$ is rejected.
	\end{assump}
	
	\begin{theorem} \label{theo:doehler_exte}
		Under arbitrary dependence we have for the SU-procedure based on the truncated critical values $\widehat \beta_{j:m}^{(k)} = \min( \widehat\alpha_{j:m},\widehat\alpha_{k:m})$ 
		\begin{align}
		\text{FDR} \leq \sum\limits_{i \in I_0} \Bigl[ \frac{P(p_i \leq \widehat{\alpha}_{k:m})}{k} +     \sum_{j=1}^{k-1} \frac{P(p_i \leq \widehat{\alpha}_{j:m})}{j(j+1)} \Bigr]. \label{eqn:theo:doehler_exte}
		\end{align}  
		The upper bound can be rewritten as
		\begin{align}
		\text{FDR} \leq \sum\limits_{i \in I_0}  \sum\limits_{j=1}^{k} \frac{P(\widehat{\alpha}_{j-1:m} < p_i \leq \widehat{\alpha}_{j:m})}{j}. \label{eqn:theo:doehler_exte_rewritten}
		\end{align}  
	\end{theorem}
	Clearly, Theorem \ref{theo:doehler_exte} includes procedures based on the original critical values $\widehat \alpha_{j:m}$, i.e., without truncation, by setting $k=m$. Without the truncation assumption \eqref{eqn:cond_hatmo} Storey's test based on \eqref{eqn:def_Storey_estimator} can be extremely liberal under strong dependence, see \citet{RomanoETAL2008}, but a sparsity modification may be helpful. Better performances can be obtained for the choice of a small tuning parameter $\lambda = \alpha/(1+\alpha)$ for Storey's procedure, see \citet{SarkarHeller2008}. Part \eqref{enu:exam_extreme_Storey} of the following Example is a modification of Proposition 17 from \citet{BlanchardRoquain2009}.
	\begin{example}[Extreme dependence]\label{example:extreme_dependence}
		Introduce $p_1=p_2=\ldots=p_{n_0} = U$ for a uniformly distributed random variabel $U$ on $(0,1)$. We set $p_{n_0+1}=\ldots=p_n=0$ for the $p$-values of false hypotheses and define $R_m(\lambda)=m \widehat F_m(\lambda)$.
		
		\begin{enumerate}[(a)]
			\item\label{enu:exam_extreme_Storey} Consider the Storey estimator \eqref{eqn:def_Storey_estimator}. The critical values are then given by
			\begin{align*}
			\widehat \alpha_{j:m} =  \frac{ \alpha j}{\widehat m_0} \wedge \lambda = 
			\begin{cases}
			\min( \alpha j (1-\alpha),\lambda) & \text{ if }U \leq \lambda, \\
			\min\Bigl(  \frac{\alpha j (1-\alpha)}{m_0+1},\lambda \Bigr) & \text{ if }U > \lambda.
			\end{cases}
			\end{align*}
			Note that we have $V=0$ for $U>s:= \min(\alpha m (1-\alpha), \lambda)$. Thus,
			\begin{align*}
			\text{FDR}= s \E \Bigl( \frac{V}{R} \Bigl | U \leq s \Bigr) = \frac{m_0}{m} \min(\alpha m (1-\alpha), \lambda),
			\end{align*}
			which equals $\frac{m_0}{m}\lambda$ when $m$ is sufficiently large. Only for $\lambda = \alpha$ we have FDR-$\alpha$-control for all $m_0$, whereas $\lambda$ is often proposed to be close to $1/2$.
			
			\item (Sparsity modification) Let $\widetilde m_0 = \max(\widehat m_0,mC)$ for some $0<C<1$ with $Cm\geq (1-\lambda)^{-1}$. Then
			\begin{align*}
			\widehat \alpha_{j:m} = \frac{ \alpha j}{\widetilde m_0} \wedge \lambda \leq \frac{ \alpha j}{Cm} \wedge \lambda
			\end{align*}
			with `` $=$`` in the case of $U\leq \lambda$. Then $\frac{V}{R}=\frac{ m_0}{m}\mathbf{1}\{U\leq \min\{(\alpha/C),\lambda\}$ and 
			\begin{align*}
			\text{FDR} = \frac{m_0}{m} \min\Bigl( \frac{\alpha}{C},\lambda \Bigr).
			\end{align*}
			A sparsity assumption expressed by the choice of $C$ close to $1$ allows the slightly increased universal FDR bound $\alpha/C$ for each $\lambda$.
		\end{enumerate}
	\end{example}

	\begin{remark}[SD tests]\label{rem:SD}
		A popular alternative to SU tests are step-down (SD) procedures. In case of the latter, the number $R$ of rejections is given by
		\begin{align*}
		R = \max\{j: p_{i:m} \leq \widehat \alpha_{i:m} \text{ for all }i\leq j\}.
		\end{align*}
		Regarding (11)--(13) in \citet{BenditkisETAL2018}, we can adjust our proofs for certain SD tests. To be more specific, the upper bound of Theorem \ref{theo:doehler_exte} is valid for any SD test based on (data driven) critical values $\eqref{eqn:hatalpha_<...}$ such that Assumption \ref{ass:constant_alpha} is fulfilled. In particular, all other FDR upper bounds from the previous sections hold. We want to point out that Assumption \ref{ass:constant_alpha} depends on the choice of $R$ and, thus, it may differ for the SU and the SD test, even when the same critical values are used.
	\end{remark}
	\section{Real data example}\label{sec:real_data}
\begin{figure}[ht]
		\begin{minipage}{0.495\textwidth}
			\begin{center}
				\includegraphics[width=1.05\textwidth, trim = 0mm 4mm 0mm 20mm,clip]{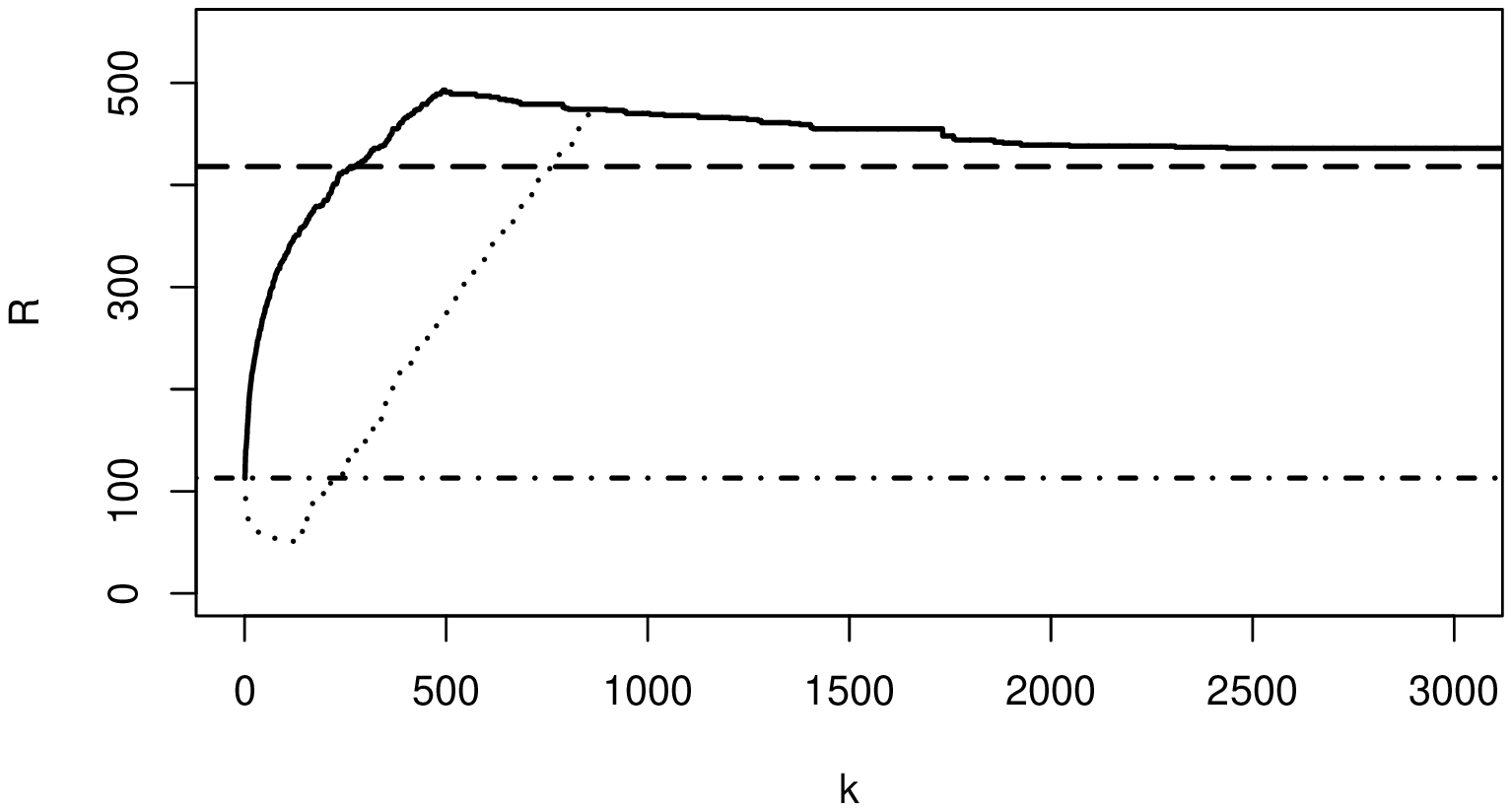}
			\end{center}
		\end{minipage}
		\begin{minipage}{0.495\textwidth}
			\begin{center}
				\includegraphics[width=1.05\textwidth,trim = 0mm 4mm 0mm 20mm,clip]{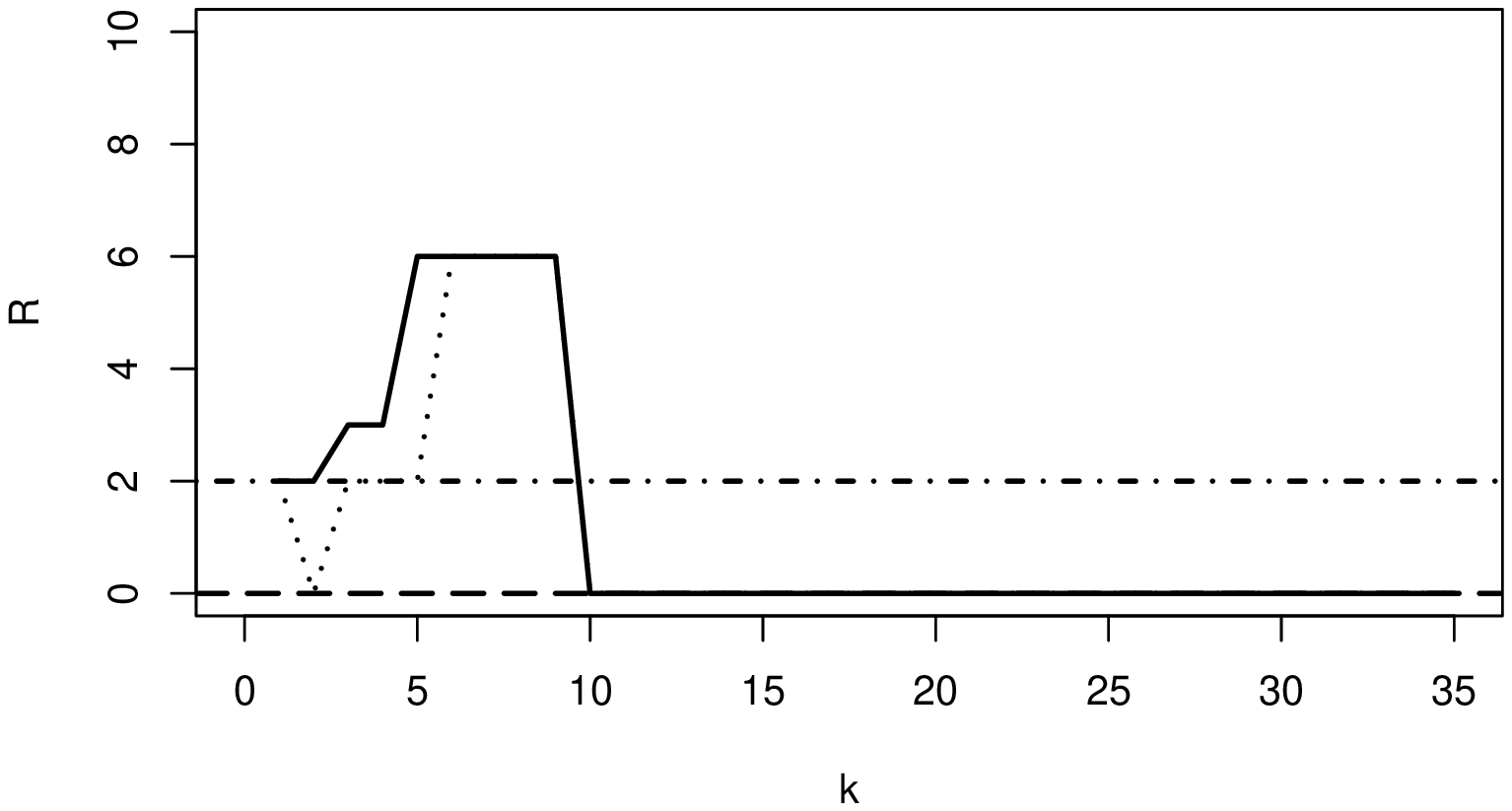}
			\end{center}
		\end{minipage}
		\caption{The number $R$ of rejections corresponding to the BH$(k)$ (solid line) and the ES$(k)$ (dotted line) procedures applying to Example \ref{exam:notterman} (left) and \ref{exam:needleman} (right) are plotted for various $k$. The height of the horizontal dashed and dot-dashed lines equal $R$ for the BY and the Bonferroni procedure, respectively. }
		\label{fig:data_example_BY}
	\end{figure}

\begin{figure}[ht]
		\begin{minipage}{0.495\textwidth}
			\begin{center}
				\includegraphics[width=1.05\textwidth, trim = 0mm 4mm 0mm 20mm,clip]{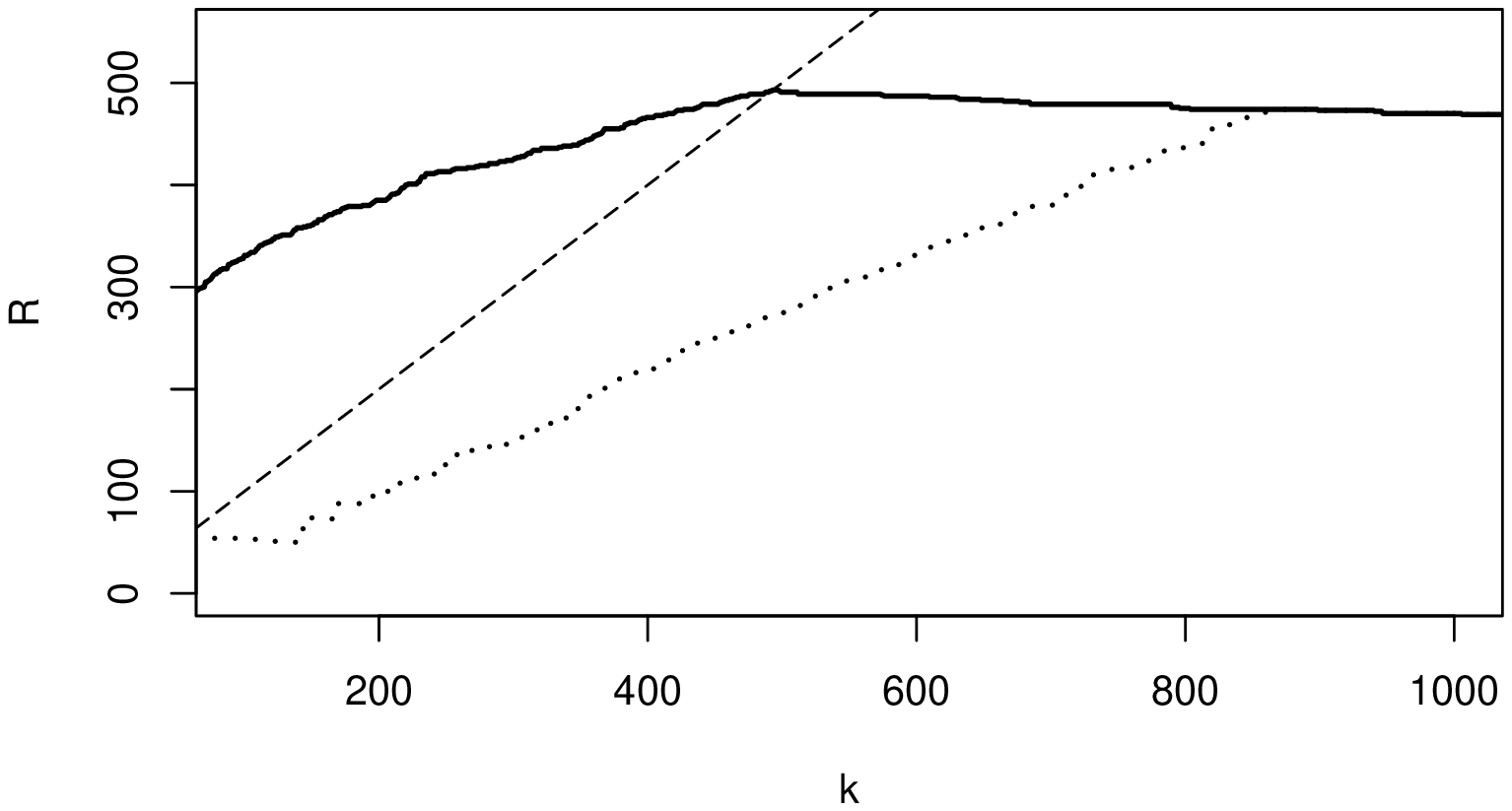}
			\end{center}
		\end{minipage}
		\begin{minipage}{0.495\textwidth}
			\begin{center}
				\includegraphics[width=1.05\textwidth,trim = 0mm 4mm 0mm 20mm,clip]{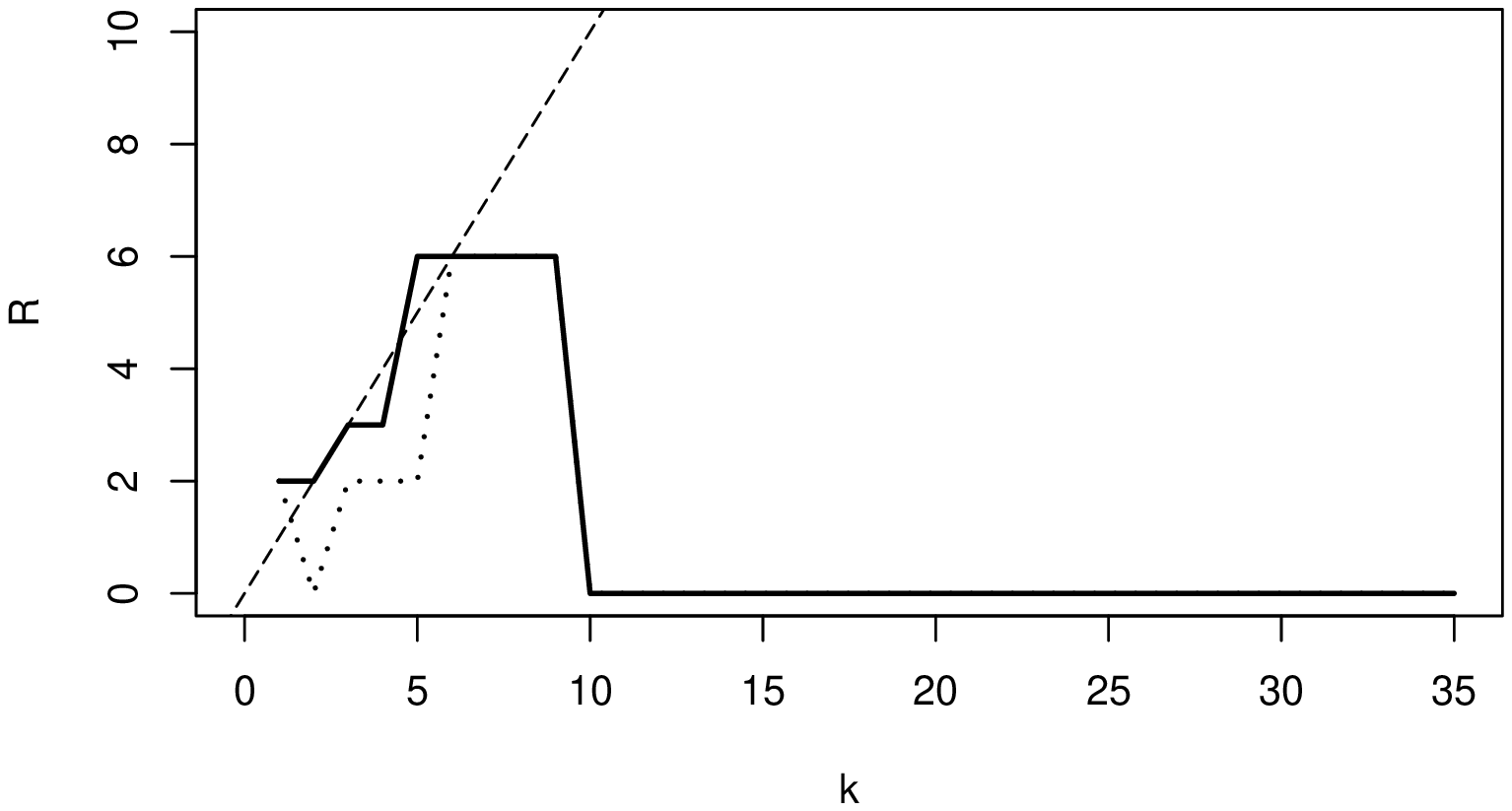}
			\end{center}
		\end{minipage}
		\caption{The number $R$ of rejections corresponding to the BH$(k)$ (solid line) and the ES$(k)$ (dotted line) procedures applying to Example \ref{exam:notterman} (left) and \ref{exam:needleman} (right) are plotted for various $k$. For comparison, the map $k\mapsto k$ is displayed (dashed line).}
		\label{fig:data_example_early_stop}
	\end{figure}

	To illustrate the benefits of our new procedures, we apply them to two different real data examples. We focus on the truncated version \eqref{eqn:def_beta} of the BH procedure, in short denoted by BH$(k)$, and the early stop procedure, in short ES$(k)$, from Section \ref{sec:early_stopp}.  These procedures are compared with the dependence corrected BH procedure (BY) of \citet{BenjaminiYekutieli2001} and the Bonferroni (Bonf) procedure, where we set $\alpha= 5\%$ for all of them.
	\begin{example}[Colon tissue]\label{exam:notterman}
		The colon adenocarcinomas data of \citet{NottermanETAL2001} consists of $7457$ gene measurements for $18$ patients on adenocarcinomas tumor and normal tissues, respectively. The data set is available, e.g., in the \textsc{R} package \textit{mutoss}. For each gene measurement we applied Welch's paired t-test resulting in $m=7457$ $p$-values.
	\end{example}
	\begin{example}[Dentine lead]\label{exam:needleman}
		\citet{needleman:etal:1979} examined the neuropsychologic effects of unidentified childhood exposure to lead. Therefore, they compared the Wechsler Intelligence Scale for Children (Revised), certain verbal processing scores, the reaction times and the classroom performances of $58$ children with high dentine lead levels and $100$ children with low dentine levels. At all, there are $m=35$ different measurements and $m$ corresponding $p$-values, see Tables 3, 7 and 8 in \cite{needleman:etal:1979}.
	\end{example}
	Since the BH procedure is only FDR-$\alpha$-controlling under specific dependence structures we exclude it in our comparison. For the sake of completeness, we nevertheless want to state that $R^{\text{BH}}=1157$ and  $R^{\text{BH}}=9$ are the number of rejections in the situation of Example \ref{exam:notterman} and Example \ref{exam:needleman}, respectively. Our truncated BH$(k)$ procedure is valid under arbitrary dependence and leads, for a variety of different $k$, to significantly more rejections than the BY and Bonferroni procedure in both examples, see Figure \ref{fig:data_example_BY}. The plot corresponding to Example \ref{exam:needleman} illustrates two general phenomenons, which should be kept in mind when choosing $k$. First, for large $k$ the improvement in comparison to the BY procedure becomes more and more negligible, which is not surprising due to the logarithmic rate of the correction factor $\sum_{i=1}^ki^{-1}$. Second, a too small $k$ may lead to fewer rejections than the BY procedure. Choosing $k$ as the truncation parameter implies that we give a special attention to the  $p$-values below $\beta_{k:m}$, see \eqref{eqn:def_beta}, ignoring the remaining ones. Hence, it is not surprising that a choice $k$ smaller than $R^{\text{BY}}$, the number of rejections for the BY procedure, may lead to fewer rejections than $R^{\text{BY}}$. In addition to expected sparse signals, the choice of $k$ can also be influenced by limited resources for further data processing steps. Under these circumstances Example \ref{exam:notterman} advocates to use the ES($k$)-procedure under restriction priorities for $R$. After some small $k$ effects the number of rejections (dotted line) is approximately linear. When the dotted line reaches the BH($k$) curve, early stopping is obsolete, see Figures \ref{fig:data_example_BY} and \ref{fig:data_example_early_stop}.
	
	%{\color{red} Du fragst, ob wir das herausstelllen sollen. M{\"o}chtest du den ganzen Absatz streichen?} We want to point out again that the BH($k$) procedure does not ensure $R\leq k$. In the situation of Example \ref{exam:needleman}, this is the case for $k=1$ (with $R=2$) and $k=5$ (with $R=6$). For Example \ref{exam:notterman} this can be observed for a variety of $k$, see Figure \ref{fig:data_example_early_stop}. But our ES$(k)$ procedure ensure $R\leq k$, at least when $k$ is larger than or equal to $R^{\text{Bonf}}$. For example, $R=2$ for the ES$(1)$ procedure in Example \ref{exam:needleman} is inline with the mentioned restriction because $R^{\text{Bonf}}=2$. 
	To judge the restriction appropriately, we want to point out that the Bonferroni test does not only control the FDR level but also the family-wise error rate $\text{FWER}=P(V>0)$, the probability that at least one rejection is false. Hence, it can be seen as a minimal requirement to beat the Bonferroni method. In this context, "to beat" means that the number of rejections should be at least $R^{\text{Bonf}}$. For small $m$, as in Example \ref{exam:needleman}, the Bonferroni method is quite powerful, e.g., it beats the BY procedure, we have $R^{\text{BY}}=0$ whereas $R^{\text{Bonf}}=2$. Although our two procedures reject up to $R=6$ hypotheses for the right choice of $k$ in Example \ref{exam:needleman}, they reject no hypotheses when $k>9$ and, in particular, they are worse than the Bonferroni method. 
	%\begin{table} 
	%	\large \centering
	%	\begin{tabular}{ c || c  c   c c c }
	%		& BH & BY & Bonf & BH$(k^*)$ & ES$(k^*)$ \\[0.1cm]
	%		\hline \hline\\[-0.3cm]
	%		Colon tissue & 1157 &418 & 113 & 493 & 493 \\[0.1cm]
	%		\hline\\[-0.3cm]
	%		Lead exposure & 9 & 0 & 2 &6
	%	\end{tabular} 
	%	\caption{The number $R$ of rejections are displayed for Example \ref{exam:notterman} (Colon tissue) and \ref{exam:needleman} (Lead exposure) when applying the classical BH, the dependence corrected BY, the Bonferroni (Bonf) and our truncated BH procedure (BH($k^*$)) with optimal truncation parameter $k^*$.}\label{tab:real_data_BY+BH}
	%\end{table}
	
	\section{Miscellaneous}\label{sec:miscellaneous}
	Under sparsity the smallest critical values are most important. For such designs we propose another promising SU-tests with a first coefficient $\alpha_{1:m}$ closer to the Bonferroni coefficient $\alpha/m$. Introduce
	\begin{align*}
	\alpha_{j:m}^{(k)} = \frac{ \alpha}{m a_k} \min \Bigl( \frac{j}{\sum_{i=1}^j i^{-1}}, \frac{k}{\sum_{i=1}^k i^{-1}} \Bigr), \quad a_k =  \Bigl( \sum_{i=1}^k \frac{1}{i} \Bigr)^{-1} + \sum_{j=1}^{k-1} \Bigl( (j+1)\sum_{i=1}^j \frac{1}{i} \Bigr)^{-1} .
	\end{align*}
	This new SU-procedure, labeled as our new sparsity test of size $k$ for short SP($k$), is motivated by Example \ref{exam:needleman}. Here, Bonferroni is superior to BY, which has too small coefficients at the beginning. Table \ref{tab:appr_alpha} displays a rough approximation of the first critical values $\alpha_{1:m}$ and $\alpha^{(k)}_{1:m}$, respectively, for different procedures when $k$ and $m$ are large.
	\begin{table} 
		\large
		\centering
		\caption{The first critical values of the new Benjamini-Yekutieli procedure as well as the ones discussed in Section \ref{sec:real_data} }\label{tab:appr_alpha}
		\begin{tabular}{ c | c ccc }
			 & Bon & BY & BH($k$) & SP($k$) \\
			\hline
			$\alpha_{1:m}$  & $= \frac{\alpha}{m}$ & $ \sim \frac{\alpha}{m\log(m)}$ & $\sim\frac{\alpha}{m\log(k)}$ & $\sim\frac{\alpha}{m\log\log(k)}$
		\end{tabular}
	\end{table}
	\begin{theorem}\label{theo:new_alpha}
		\begin{enumerate}[(a)]
			\item\label{enu:theo:new_alpha_BI} (cf. \cite{HeesenJanssen2015}, Prop. 4.1) Under BI we have
			\begin{align*}
			\text{FDR} \leq \frac{m_0}{m} \frac{ \alpha}{a_k}.
			\end{align*}
			
			\item\label{enu:theo:new_alpha_dep} Under arbitrary dependence $\text{FDR}\leq (m_0/m)\alpha$.
			
			\item\label{enu:theo:new_alpha_norm_coeff} For $k\geq 4$ the normalizing coefficients obey the inequality
			\begin{align*}
			 	\log(1+\log(k+1)) - \log( 1 + \log 3) + \frac{1}{1+ \log(k)} \leq a_k \leq \frac{13}{18} + \log\log k + \frac{1}{\log(k+1)}.
			\end{align*}
		\end{enumerate}
	\end{theorem}
	Since the degree of sparsity is typically unknown the truncation parameter $k$ of BH($k$) should be not too small. But, for larger $k$ the BH($k$) coefficients are smaller, which should be judged within a trade-off. Here, the new SP($k$) procedure helps with increased lower coefficients. The substitution of BH($k$) by SP($k$) then gives more security for very sparse signals. Note that for each $k$ there exists a largest value $j_0$ with
	\begin{align}\label{eqn:NBY_j0}
		\alpha^{(k)}_{j:m}\geq \beta_{j:m}^{(k)}\quad \text{ iff }j\leq j_0.
	\end{align}	
	 For intermediate $k$ the crucial index $j_0$ can be obtained directly. Using the approximations $\sum_{i=1}^ki^{-1} \sim \log(k)$ and $a_k\sim \log\log(k)$ we get approximately $j_0\sim k^{1/\log\log(k)}$. In contrast to BH($k$), the SP($k$) procedure is at least as good as the Bonferroni test for every $k$ in the context of Example \ref{exam:needleman}.
	
	\section{Proofs}\label{sec:proofs}
	We first present the proofs corresponding to our main results from Section \ref{sec:main}. After that we prove the remaining statements in order of their appearance in the paper. 
	\subsection{Proof of Theorem \ref{theo:main_res}}
	\underline{\eqref{enu:theo:main_res_BI}:} First, observe that \citet{HeesenJanssen2015} already proved, see their Lemma 7.1(a), that
	\begin{align*}
	\E \Bigl( \frac{V}{m\widehat \alpha_{R:m}} \Bigr)	= \frac{m_0}{m}.
	\end{align*}
	By our assumptions we have $g(g^{-1}(\lambda))\leq \lambda$ and
	\begin{align*}
	\widehat \alpha_{R:m} = \min \Bigl( g\Bigl( \frac{R\wedge k}{\widehat m_0}\Bigr),\lambda  \Bigr) \leq g\Bigl( \frac{R}{\widehat m_0}\wedge B\Bigr)\quad 
	\end{align*}
	For $R\geq 1$ we can conclude
	\begin{align*}
	\frac{ R }{ \widehat m_0 \widehat \alpha_{R:m}}\geq \frac{ B \wedge (R/\widehat m_0) }{ g(B \wedge (R/\widehat m_0))} \geq  \Bigl( \sup_{t\in I} \frac{ g(t) }{ t } \Bigr)^{-1}.
	\end{align*}
	Finally, we can deduce
	\begin{align*}
	\E \Bigl(  \frac{ V}{R} \Bigr) \Bigl( \sup_{t\in I} \frac{ g(t) }{ t }\Bigr)^{-1} &\leq \E \Bigl( \frac{V}{\widehat m_0 \widehat \alpha_{R:m}} \Bigr) \leq \frac{1}{C} \E \Bigl( \frac{V}{m \widehat \alpha_{R:m}} \Bigr) = \frac{m_0}{mC},
	\end{align*}
	which completes the proof of \eqref{eqn:main_res_BI}. Analogously, the statement for the deterministic case can be obtained.\\

	The basic idea of the proof of \eqref{enu:theo:main_res_dep_det} and \eqref{enu:theo:main_res_dep_adap} is to rewrite  the generating function $g$ as
	\begin{align}\label{eqn:g=int_mu}
	g\Bigl( \frac{x}{m} \Bigr) = \frac{1}{m} \int_0^x z \,\mathrm{ d }\mu(z)
	\end{align}
	for $\mu$-almost all $x\in[0,m]$ and a finite measure $\mu$ specified later, where $\mu$ may depend on $m$. It should be mentioned that \citet{BlanchardRoquain2009} already used this technique. 
	
	\underline{\eqref{enu:theo:main_res_dep_det}}: For the generating function $g$ the measure $\mu$ needs to be specified only on $D=[1,mB]$ if $\mu((0,1/m)\cup(mB,m))=0$. We denote by $\epsilon_x$ the Dirac measure centred at $x$, i.e., $\epsilon_x(A)=\mathbf{1}\{x\in A\}$. Let
	\begin{align*}
	\mu = m\sum_{1\leq j \leq mB} \frac{ g(j/m) - g((j-1)/m)}{j} \epsilon_j.
	\end{align*}
	Then \eqref{eqn:g=int_mu} holds, which implies for all  $j=1,\ldots,m$ that
	\begin{align*}
	\alpha_{j:m} = g\Bigl( \frac{j}{m} \Bigr) = \frac{\mu(D)}{m} \int_1^j x \,\mathrm{ d } \nu(x)
	\end{align*}
	for the probability measure $\nu=\mu/\mu(D)$. Lemma 7.1(b) of \citet{HeesenJanssen2015} with $\alpha = \mu(D)$ and $\widehat \rho(i)=i$ yields $\text{FDR}\leq (m_0/m)\mu(D)$. Finally, note  
	\begin{align*}
	\mu(D) &= \sum_{1\leq j \leq mB} m\frac{ g(j/m) - g((j-1)/m)}{j} \\
	&= \frac{ g([mB]/m)}{[mB]/m} + \sum_{j=1}^{[mB]-1} \frac{ g(j/m)}{(j+1)j/m}  \leq \Bigl( \sum_{1\leq j \leq mB} \frac{1}{j} \Bigr) C_k^{\text{det}}.
	\end{align*}
	
	\underline{\eqref{enu:theo:main_res_dep_adap}}: Since $j/\widehat m_0 \geq \delta/m$ we just need to specify $\mu$ on $D=[\delta,mB]$. Here, we define $\mu$ as a mixture of a Dirac measure and a Lebesgue measure with density $h$:
	\begin{align*}
	\mu = \frac{m}{\delta} g\Bigl( \frac{ \delta }{m} \Bigr)\epsilon_\delta + (h \mathbf{1}_{[\delta,mB]}) \cdot \lebesgue \quad\text{with } h(z) = \frac{1}{z} g'\Bigl( \frac{z}{m} \Bigr).
	\end{align*}
	Setting again $\nu=\mu/\mu(D)$ we obtain
	\begin{align*}
	&\widehat \alpha_{j:m}=g\Bigl( \frac{j}{\widehat m_0} \Bigr) = \frac{\mu(D)}{m} \int_\delta^{j\frac{m}{\widehat m_0}} z \,\mathrm{ d }\nu(z).
	\end{align*}
	Applying Lemma 7.1(b) of \citet{HeesenJanssen2015}, now with $\widehat \rho(j)= j (m / \widehat m_0)$, yields 
	\begin{align*}
	\frac{ \text{FDR} }{ 1/C}\leq \E\Bigl( \frac{V}{ (m/\widehat m_0)R} \Bigr) \leq \frac{m_0}{m} \mu(D).
	\end{align*}
	Finally, using integration by parts we can deduce
	\begin{align*}
	\mu(D) & = \frac{m}{\delta} g\Bigl( \frac{ \delta}{m} \Bigr) + \int_{\delta}^{mB} \frac{1}{z}g'\Bigl( \frac{z}{m} \Bigr)\,\mathrm{ d }z  = \frac{g(B)}{B}  + \int_\delta^{mB}\frac{m}{z^2}g\Bigl( \frac{z}{m} \Bigr)\,\mathrm{ d }z \\
	&\leq \sup_{(\delta/m)\leq t \leq B} \frac{g(t)}{t} \Bigl[ 1 + \int_\delta^{mB} \frac{1}{z} \,\mathrm{ d }z  \Bigr].
	\end{align*}

	\subsection{Proof of Theorem \ref{theo:doehler_exte}} 
	To simplify the notation, assume  $\widehat \alpha_{j:m} = \widehat \beta_{j:m}^{(k)}$. Let $p=(p_1,\ldots,p_m)$ and set $p^{(i)} = (p_1,\ldots,p_{i-1},0,p_{i+1},\ldots,p_m)$, where only the $i$-th coordinate $p_i$ is replaced by $0$ for $i \in I_0$. For the proof, we prefer to explicitly state the dependence from $p$ and write $\widehat{\alpha}_{j:m}(p),V(p),R(p)$. Assumption \ref{ass:constant_alpha} implies for all $j$
	\begin{align}\label{eqn:proof_doehl_alpha}
	\widehat{\alpha}_{j:m}(p) = \widehat{\alpha}_{j:m}(p^{(i)}) \textrm{ on the set } \{p_i \leq \widehat{\alpha}_{R(p):m}(p)\}.
	\end{align} 
	Introduce for fixed $i\in I_0$ the sets $D_{0}^{(i)}(p^{(i)}) = \emptyset$ and \[D_{r}^{(i)}(p^{(i)})=\{R(p^{(i)})=r \} \textrm{ for } 1 \leq r \leq m. \] Below, we use the well-known inclusion for step-up tests:
	\begin{align}\label{eqn:proof_doehl_sets}
	\{ p_i \leq \widehat{\alpha}_{r:m}(p) \textrm{ and } R(p)=r \} \subset \{ p_i \leq \widehat{\alpha}_{r:m}(p) \} \cap D_{r}^{(i)}(p^{(i)}).
	\end{align}
	For $i\in I_0$ we have
	\begin{align*}%\label{eqn:proof_doehl_Ei}
	E_i := E \! \left( \frac{\mathbf{1}(\textrm{"$H_i$ is rejected"})}{R(p)} \right) = \sum\limits_{r=1}^{m} \frac{1}{r} P\! \left(\left\{ p_i \leq \widehat{\alpha}_{r:m}(p) \right\} \cap \left\{ R(p) =r \right\} \right). 
	\end{align*}
	The $r$-th summand can be transformed by \eqref{eqn:proof_doehl_sets}:
	\begin{align*}
	& P  \Bigl( \Bigl\{ p_i \leq \widehat{\alpha}_{r:m}(p), R(p) = r \Bigr\} \Bigr) \leq P \Bigl( \Bigl\{ p_i \leq \widehat{\alpha}_{r:m}(p) \Bigr\} \cap D_{r}^{(i)}(p^{(i)}) \Bigr) \nonumber\\
	&= \sum\limits_{j=1}^{r} P  \Bigl( \Bigl\{ \widehat{\alpha}_{j-1:m}(p) < p_i \leq \widehat{\alpha}_{j:m}(p)\Bigr\} \cap D_{r}^{(i)}(p^{(i)}) \Bigr) \nonumber\\
	&=: \sum\limits_{j=1}^{r} p_{irj}. \nonumber
	\end{align*}
	Interchanging the summation yields
	\begin{align}
	E_i &= \sum\limits_{r=1}^{m} \sum\limits_{j=1}^{r} \frac{1}{r} p_{irj} = \sum\limits_{j=1}^{m} \sum\limits_{r=j}^{m} \frac{1}{r} p_{irj} \nonumber\\
	&\leq \sum\limits_{j=1}^{m} \sum\limits_{r=j}^{m} \frac{1}{j} p_{irj} \leq \sum\limits_{j=1}^{m} \frac{1}{j} \sum\limits_{r=j}^{m}  p_{irj} \nonumber\\
	&\leq \sum\limits_{j=1}^{m} \frac{1}{j} P \! \left(\widehat{\alpha}_{j-1:m} < p_i \leq \widehat{\alpha}_{j:m} \right). \nonumber
	\end{align}
	The equality $\text{FDR} = \sum_{i \in I_0} E_i$ proves \eqref{eqn:theo:doehler_exte_rewritten} because $\widehat \alpha_{j:m}=\widehat \alpha_{k:m}$ for $j>k$. The other representation of the upper bound, see \eqref{eqn:theo:doehler_exte}, follows from a rearrangement and the trivial observation
	\begin{align*}
	P \left( \widehat{\alpha}_{j-1:m} < p_i \leq \widehat{\alpha}_{j:m} \right) = P\left(p_i \leq \widehat{\alpha}_{j:m} \right) - P \left(p_{i} \leq \widehat{\alpha}_{j-1:m}\right).
	\end{align*} 
	
	\subsection{Proof of Theorem \ref{theo:mod_Blan+Roq}}
	The proof is related to \cite{HeesenJanssen2016} and relies on conditioning with respect to the $\sigma$-algebra
	\begin{align*}
	\mathcal{F}_{\lambda} := \sigma\Bigl( \mathbf{1}\{p_{i}\leq s\}: s\in[\lambda,1], \, 1 \leq i \leq m \Bigr).
	\end{align*}
	Adopting the notation of the proof for Theorem \ref{theo:doehler_exte} we obtain from Fubini's Theorem that for $i\in I_0$
	\begin{align} \label{eqn:proof_prop}
	&\E \Bigl( \frac{ \mathbf{1}\{p_i \leq \alpha_{R:m}\}}{R} \Bigl | \mathcal F_\lambda \Bigr) \mathbf{1}\{p_i\leq \lambda\} \\
	&=\E \Bigl( \frac{ \mathbf{1}\{p_i \leq \alpha_{R(p^{(i)}):m}\}}{R(p^{(i)})} \Bigl | \mathcal F_\lambda \Bigr) \mathbf{1}\{p_i\leq \lambda\} \nonumber \\ &=\E \Bigl( \frac{ \alpha_{R(p^{(i)}):m}}{R(p^{(i)})} \Bigl | \mathcal F_\lambda \Bigr) \mathbf{1}\{p_i\leq \lambda\} \nonumber \\
	&\leq (1-\lambda)\alpha\E \Bigl( \frac{ 1}{m-R(p^{(i)})(1-\alpha)} \Bigl | \mathcal F_\lambda \Bigr) \mathbf{1}\{p_i\leq \lambda\}, \nonumber 
	\end{align}
	see also \cite{BenditkisETAL2018}. First, consider the case $m_0=1$. Taking the expectation of \eqref{eqn:proof_prop} we obtain from the trivial fact $R(p^{(i)}) \leq m$ that
	\begin{align*}
	\text{FDR} \leq   \frac{ (1-\lambda)\alpha \lambda}{m - R(p^{(i)})(1-\alpha)}  \leq \frac{(1-\lambda)\alpha \lambda}{m \alpha} \leq \alpha.
	\end{align*}
		
	Under $m_0\geq 2$ the random variable $R(p^{(i)})$ is bounded by $R(p^{(i)}) \leq m_1+1+V^{(i)}(\lambda)$ with $V^{(i)}(\lambda)= \# \{j\in I_0\setminus \{i\}: p_j\leq \lambda\}$, where $m_1$ is the number of false hypotheses.  Since
	\begin{align*}%\label{eqn:proof_prop2}
	m - R(p^{(i)})(1-\alpha) &\geq \alpha m + (1-\alpha)[(m_0-1)-V^{(i)}(\lambda)] \nonumber \\
	& = (1-\alpha) \Bigl[ K + (m_0-1) - V^{(i)}(\lambda) \Bigr]\quad \text{for }K:= \frac{\alpha m}{1-\alpha}
	\end{align*}
	we can bound \eqref{eqn:proof_prop} by
	\begin{align*}
	\frac{(1-\lambda)\alpha}{ 1 - \alpha } \E \Bigl( \frac{ 1 }{ K + X^{(i)} } \Bigl | \mathcal F_\lambda \Bigr) \mathbf{1}\{p_i\leq \lambda\},\quad \text{where }X^{(i)}:=(m_0-1)- V^{(i)}(\lambda).
	\end{align*}
	Since $X^{(i)}$ and $p_i$ are independent for all $i\in I_0$, taking expectations and summing up over $i\in I_0$ we obtain
	\begin{align}\label{eqn:proof_prop_FDR}
	\text{FDR} \leq \frac{ (1-\lambda)\alpha\lambda}{1-\alpha} \sum_{i\in I_0}\E\Bigl( \frac{1}{K + X^{(i)}} \Bigr).
	\end{align}
	Clearly, $V^{(i)}(\lambda)$ is binomial distributed with parameters $m_0-1$ and $\lambda$ under BI. Thus, we can use the following well-known inequalities
	\begin{align*}
	\E\Bigl( \frac{ V^{(i)}(\lambda)}{ 1 + X^{(i)}} \Bigr) \leq \frac{ \lambda}{ 1 - \lambda}\quad\text{and}\quad \E\Bigl( \frac{1}{1+X^{(i)}} \Bigr) \leq \frac{1}{ m_0(1-\lambda)}.
	\end{align*}
	It is easy to check that 
	\begin{align*}
	t\mapsto h(t) := \frac{ t }{ 1+ t(K-1)}
	\end{align*}
	is increasing and concave with $h((t+1)^{-1})=(K+t)^{-1}$. Thus, Jensen's inequality implies
	\begin{align}\label{eqn:proof_prop_EK+X}
	\E\Bigl( \frac{1}{K+X^{(i)}} \Bigr) \leq h\Bigl( \E \Bigl( \frac{1}{1+X^{(i)}} \Bigr) \Bigr) \leq h \Bigl( \frac{1}{m_0(1-\lambda)} \Bigr).
	\end{align}
	Combining \eqref{eqn:proof_prop_FDR} and \eqref{eqn:proof_prop_EK+X} with $K-1=\frac{(m+1)\alpha -1}{1-\alpha}$ yields the FDR bound
	\begin{align*}
	\text{FDR} \leq \alpha \lambda \Bigl( 1- \alpha + \frac{ (m+1)\alpha - 1 }{ m_0(1-\lambda)} \Bigr)^{-1}.
	\end{align*}
	Since $m_0\leq m$ it is sufficient to prove 
	\begin{align}\label{eqn:proof_leq_lam}
	(1-\alpha) + \frac{(m+1)\alpha - 1}{ m_0(1-\lambda)} \geq \lambda.
	\end{align}
	By straight forward calculations, \eqref{eqn:proof_leq_lam} holds if and only if
	\begin{align}\label{eqn:proof_iff}
	(1-\lambda)^2m - \alpha m (1-\lambda) + (m+1)\alpha - 1 \geq 0.
	\end{align}
	To verify this inequality, let us introduce the quadratic function $\psi:\R\to\R$ given by $\psi(x):=mx^2+(m+1)\alpha - 1 -\alpha m x \geq 0$. Obviously, this function attains its minimum at $x=\alpha/2$ with
	\begin{align*}
	\psi\Bigl( \frac{\alpha}{2} \Bigr) = m  \Bigl( \alpha - \frac{ \alpha^2}{4} \Bigr) + \alpha - 1 \geq 0.
	\end{align*}
	This proves \eqref{eqn:proof_iff} and, finally, completes the proof.

	\subsection{Proof of Theorem \ref{theo:early_stopp}}
	\underline{\eqref{enu:theo:early_stopp_BI}}: Let $i\in I_0$. We adopt the notation of the proof of Theorems \ref{theo:mod_Blan+Roq} and \ref{theo:doehler_exte}.  As already noted there, we have $\widehat \alpha_{j:m}(p) = \widehat \alpha_{j:m}(p^{(i)})$ and $R(p)=R(p^{(i)})$ whenever $p_i\in[0,\widehat \alpha_{R(p):m}]$. Thus, we can deduce from Fubini's Theorem that
	\begin{align*}
	\text{FDR} &= \sum_{i\in I_0} \E \Bigl( \frac{ \mathbf{1}\{p_i \leq \widehat \alpha_{R(p):m}(p)\}}{R(p)} \Bigr) = \sum_{i\in I_0} \E \Bigl( \frac{ \mathbf{1}\{p_i \leq \widehat \alpha_{R(p^{(i)}):m}(p^{(i)})\}}{R(p^{(i)})} \Bigr) \\
	&= \sum_{i\in I_0} \E \Bigl( \frac{ \widehat \alpha_{R(p^{(i)}):m}(p^{(i)})}{R(p^{(i)})} \Bigr) \leq \frac{m_0}{m} \alpha.
	\end{align*}
	
	\underline{\eqref{enu:theo:early_stopp_dep}}: Combining $\widehat\alpha_{j:m} \leq (\alpha/m)\min(j,k)$ and Theorem \ref{theo:doehler_exte} verifies the statement.

	\subsection{Proof of Theorem \ref{theo:new_alpha}}
	\underline{\eqref{enu:theo:new_alpha_BI}:} The statement follows immediately from Proposition 4.1 of \citet{HeesenJanssen2015}.
	
	\underline{\eqref{enu:theo:new_alpha_dep}:} By Theorem \ref{theo:doehler_exte}
	\begin{align*}
	\text{FDR} \leq \frac{m_0\alpha}{a_km} \Bigl( \frac{1}{\sum_{i=1}^k i^{-1} } + \sum_{j=1}^{k-1} \frac{1}{(j+1)\sum_{i=1}^j i^{-1} } \Bigr).
	\end{align*}
	
	\underline{\eqref{enu:theo:new_alpha_norm_coeff}:} From $\log(1+j)\leq \sum_{i=1}^j i^{-1} \leq 1 + \log j$ we obtain 
	\begin{align*}
	\sum_{j=3}^{k-1} \frac{1}{(j+1)\sum_{i=1}^j i^{-1} } &\leq \sum_{j=3}^{k-1} \frac{1}{ (j+1) \log(j+1)}
	\leq \int_3^k \frac{1}{x \log x} \,\mathrm{ d }x \\
	&= \log\log k - \log \log 3 \leq \log\log k
	\end{align*}
	and, consequently, the upper bound for $a_k$. Similar arguments lead to the lower bound:
	\begin{align*}
	\sum_{j=1}^{k-1} \frac{1}{(j+1)\sum_{i=1}^j i^{-1} } &\geq \sum_{j=1}^{k-1} \frac{1}{(j+1)(1+\log(1+j)) }\\
	&\geq \int_2^k \frac{1}{(x+1)(1+\log(x+1))} \,\mathrm{ d } x \\
	&= \log( 1 + \log(k+1)) - \log(1+\log 3).
	\end{align*}
	
	\bibliographystyle{elsarticle-harv} 
	\bibliography{sample}
	
\end{document}